\documentclass{conm-p-l}

\usepackage{amssymb, latexsym, amsmath, amscd, amsthm, verbatim, graphicx, xcolor, hyperref, mathtools,booktabs}

\theoremstyle{definition}
\newtheorem{remark}{Remark}

\definecolor{mylinkcolor}{rgb}{0.65,0.0,0.0}
\definecolor{myurlcolor}{rgb}{0.0,0.0,0.65}
\hypersetup{colorlinks=true,urlcolor=myurlcolor,citecolor=myurlcolor,linkcolor=mylinkcolor,linktoc=page,breaklinks=true}

\newcommand{\Z}{{\mathbb{Z}}}
\newcommand{\Q}{{\mathbb{Q}}}
\DeclareMathOperator{\Tr}{Tr}
\DeclareMathOperator{\GL}{GL}
\DeclareMathOperator{\Jac}{Jac}

\begin{document}

\title{Completely decomposable modular Jacobians}

\author{Jennifer Paulhus}
\address{Mount Holyoke College}
\email{jpaulhus@mtholyoke.edu}

\author{Andrew V. Sutherland}
\address{Massachusetts Institute of Technology}
\email{drew@math.mit.edu}
\thanks{Sutherland was supported by Simons Foundation grant 550033.}


\date{January 31, 2025}

\begin{abstract} We use recently developed algorithms and a new database of modular curves constructed for the {\it L-functions and Modular Forms Database} to enumerate completely decomposable modular Jacobians of level $N < 240$.  In particular, we find examples in 13 previously unknown genera of Jacobian varieties isogenous to a product of elliptic curves over $\Q$.  The new genera are: $38$, $68$, $75$, $76$, $77$, $78$, $113$, $135$, $137$, $157$, $159$, $169$, and $409$.
\end{abstract}

\maketitle

As abelian varieties, Jacobian varieties can be decomposed (up to isogeny) into a product of simple abelian varieties. Many interesting questions can be asked about these decompositions. In particular, Ekedahl and Serre \cite{ekedahl-serre} studied the question of when a Jacobian variety factors into a product of elliptic curves (abelian varieties of dimension one), called a {\it completely decomposable Jacobian}. They asked if the genera for which there is such a curve are bounded, and if not, whether there is an example in every genus of a curve whose Jacobian is completely decomposable.  Toward the latter question, they produced curves of many genera (up to genus 1297) with completely decomposable Jacobians defined over number fields.  Both questions remain open over fields of characteristic zero (in positive characteristic one can realize infinitely many genera; see \cite{char-p-1, char-p-2, char-p-3}).

But there were still many genera below 1297 for which Ekedahl and Serre were unable to produce an example.  In the intervening years some gaps have been filled in. Work of Yamauchi \cite{yamauchi} on decompositions of Jacobian varieties of modular curves found some new examples. The first author in collaboration with Rojas \cite{paulhus-rojas} used classifications of the automorphism groups of curves as well as a new technique involving intermediate quotients of curves with automorphisms to produce many more examples. Rojas and Rodr\'{i}guez recently demonstrated an example in genus 101 \cite{anitarubi}. See \cite[\S 3.11]{serre-book} for a brief summary of this problem and prior results.

In this short note we present results obtained from a new database of modular curves constructed for the L-functions and Modular Forms Database (LMFDB) \cite{lmfdb} to obtain completely decomposable Jacobians in 13 new genera, all of which are defined over $\Q$.  In particular, we find an example in genus 38, which increases the least genus for which no completely decomposable Jacobian over a number field (or even over $\mathbb C$) is known to 56.  The complete list of genera known to occur in characteristic zero now stands as follows (with the 13 new genera in bold):
\begin{center}
0--37, \textbf{38}, 39--55, 57--58, 61--65, 67, \textbf{68}, 69, 71--73, \textbf{75--78}, 79--82, 85, 89, 91, \\ 
93, 95, 97, 101, 103, 105--107, 109, \textbf{113}, 118, 121, 125, 129, \textbf{135}, \textbf{137}, 142,  \\
145, 154, \textbf{157}, \textbf{159}, 161, 163, \textbf{169}, 193, 199, 205, 211, 213, 217, 244, 257, \\  325, \textbf{409}, 433, 649, 1297.
\end{center}

Many of the examples produced by Ekedahl and Serre in \cite{ekedahl-serre} are Jacobians $J_0(N)$ of classical modular curves $X_0(N)$ whose non-cuspidal points parameterize elliptic curves with a cyclic isogeny of degree $N$.  They noted that $X_0(N)$ is completely decomposable whenever the eigenvalues of the Hecke operators $T_p$ acting on weight-2 cusp forms of level $N$ are rational for all $p\nmid N$, and used an early precursor to the LMFDB (unpublished tables of modular forms computed by Cohen, Skoruppa, and Zagier in the early 1990s), to determine the $N\le 1000$ satisfying these conditions. Yamauchi later determined the entire list of $N$ for which $J_0(N)$ is completely decomposable \cite[Thm 1.1]{yamauchi}, modulo errors at $N=28$ (inadvertently omitted) and $N=672$ (incorrectly included), as noted in \cite[p. 82]{EHR}. The genera that arise are:
\begin{center}
0--11, 13, 17, 19, 21, 25, 29, 33, 37, 43, 49, 53, 55, 57, 61, 73, 97, 121, 161, and 205.
\end{center}

Recently, a database of modular curves was added to the beta version of the LMFDB (\url{https://beta.lmfdb.org/ModularCurve/Q/}) \cite{lmfdb} as part of a project supported by the Simons Collaboration in Arithmetic Geometry, Number Theory, and Computation.
This database includes all modular curves $X_H$ associated with open $H\le \GL_2(\widehat\Z)$ of level $N\le 70$ with $\det(H)=\widehat \Z^\times$.  The \textit{level} of $H$ is the least integer $N$ for which $H$ is the inverse image of its projection to $\GL_2(\Z/N\Z)$, and the constraint on $\det(H)$ ensures that $X_H/\Q$ is a \textit{nice} (smooth, projective, geometrically integral) curve.  The $X_H$ are generalizations of the classical modular curves $X_0(N)$ and $X_1(N)$, and they  parameterize elliptic curves with $H$-\textit{level structure}; we refer the reader to \cite[\S 2]{RSZB} for further details.

This database includes the decomposition of $J_H\coloneqq \Jac(X_H)$ into simple modular abelian varieties $A_f$ associated to Galois orbits of weight-2 newforms $f$.
As proved in \cite{RSZB}, all such $f$ are (not necessarily new) normalized eigenforms in $S_2(\Gamma_1(N)\cap\Gamma_0(N^2))$.
These decompositions were computed using the algorithm described in \cite[\S 6]{RSZB}, which we briefly recall.  If $[f_1],\ldots,[f_m]$ is a complete list of Galois orbits of eigenforms in $S_2(\Gamma_1(N)\cap\Gamma_0(N^2))$, then there is a unique sequence of nonnegative integers $e_1,\ldots,e_m$ for which
\begin{equation}\label{eq:LA}
L(J_H,s) = \prod_{i=1}^m L(A_{f_i},s)^{e_i},
\end{equation}
which implies the integer linear constraint $a_p(J_H) = \sum_{i=1} e_i a_p(\Tr(f_i))$ for each prime~$p$.  Here $\Tr(f_i)$ denotes the \emph{traceform} of $f_i$, the sum of the Galois conjugates of $f$, see \cite[\S 4.5]{computing_cmf}).  We also have the constraint $\sum_{i=1}^m e_ia_1(\Tr(f_i)) = g(X_H)$.

The integer $q$-expansions $\sum_{n=1} a_n q^n$ of the traceforms $\Tr(f)$ were computed using the algorithms described in \cite{computing_cmf} for $n\le B$.  By increasing the bound~$B$ as required, one eventually obtains a linear system determined by the constraints on $a_p(J_H)$ and  $g(X_H)$ that has a unique solution in $\Q$ (which is necessarily integral).  This reduces the problem of computing the isogeny decomposition of $J_H$ to linear algebra, provided one can efficiently compute the integers $a_p(J_H)$ and the $q$-expansions of all the relevant traceforms $Tr(f_i)$.  An efficient algorithm for computing $a_p(J_H)$ is described in \cite[\S 5]{RSZB}, and a Magma implementation of this algorithm is available in the associated GitHub repository.  The computations of the traceforms is more difficult (this requires decomposing spaces of modular forms of level $N^2$, including spaces with nontrivial character), but for $H$ of level $N\le 70$ most of the required traceforms were already present in the LMFDB thanks to \cite{computing_cmf} (this work involved hundreds of CPU years of computations that we fortunately did not need to repeat).  This made it a simple matter to search the database for examples where $J_H$ is completely decomposable, leading to the discovery of 11 previously unknown genera.

The database also includes many $X_H$ with $H$ of level $N>70$, but for $N>70$ it does not include decompositions of the Jacobians $J_H$ due to the difficulty of computing the modular forms in $S_2(\Gamma_1(N) \cap \Gamma_0(N^2))$. Indeed, $N=71$ is already problematic in this regard, as the space of newforms with LMFDB label \href{https://www.lmfdb.org/ModularForm/GL2/Q/holomorphic/5041/2/g/}{\texttt{5041.2.g}} of dimension 9120 has resisted all efforts to decompose it to date.  However, the modular abelian varieties $A_f$ corresponding to Galois orbits of newforms $f$ in this space all have dimension divisible by 24 (the degree of the field of values $\Q(\zeta_{35})$ of the Dirichlet character orbit \href{https://www.lmfdb.org/Character/Dirichlet/5041/g}{\texttt{5041.g}}), which implies that they cannot appear as isogeny factors of a completely decomposable $J_H$.  Only weight-2 newforms with rational coefficients, corresponding to elliptic curves $E/\Q$, are relevant to our search, and these $E$ are easy to determine.  Indeed, the LMFDB contains all elliptic curves $E/\Q$ of conductor less than 500,000, which is enough to handle any $N\le 707$.

It is not computationally feasible to enumerate all open $H\le \GL_2(\widehat\Z)$ of level $N\le 707$; at level 80 there are already nearly 20 million $H$ with $\det(H)=\widehat\Z^\times$.  But it is not necessary to enumerate all such $H$.  If $J_{H}$ is not completely decomposable, then neither is $J_{K}$ for any $K\le H\le \GL_2(\widehat\Z)$, since the inclusion $K\le H$ induces a morphism of curves $X_K\to X_H$ that makes $J_H$ an isogeny factor of $J_K$.  
This suggests the following strategy: for a given level $N$, enumerate the subgroup lattice of $\GL_2(\Z/N\Z)$ up to conjugacy from the top down by successively computing (conjugacy classes of) maximal subgroups (a task that \textsc{Magma} \cite{magma} performs very efficiently), and do not proceed beyond any $H$ for which $J_H$ is not completely decomposable.
An open subgroup $H\le \GL_2(\widehat\Z)=\varprojlim_N \GL_2(\Z/N\Z)$ of level $N$ is  uniquely determined by its projection to $\GL_2(\Z/N\Z)$, and we may identify $H$ with the conjugacy class of its image in $\GL_2(\Z/N\Z)$, since conjugate subgroups give rise to isomorphic modular curves. Not every subgroup $K\le \GL_2(\Z/N\Z)$ has level $N$, it may have level $M$ properly dividing $N$, but such $K$ are still of interest because they may contain subgroups $H$ of level $N$; as noted above, if $J_K$ is not completely decomposable then neither is $J_H$.

As we enumerate the subgroup lattice of $\GL_2(\Z/N\Z)$ we may restrict our attention to $H$ for which $\det(H)=(\Z/N\Z)^\times$, as we are interested in nice curves over $\Q$ (if $\det(H)\ne (\Z/N\Z)^\times$ then $X_H$ will not be geometrically irreducible), and we can also require $H$ to contain $-I$, since if $H$ does not contain $-I$ and we put $H'=\langle H,-I\rangle$, the curves $X_H$ and $X_H'$ will be isomorphic (even though they have different moduli interpretations); see \cite[\S 2]{RSZB} for further details.  Note that if $H$ does not satisfy both these constraints than neither does any of its subgroups.

When $J_H$ is not completely decomposable we can prove this using the $L$-series of the elliptic curves $E/\Q$ of conductor dividing $N^2$.  Explicitly, let $E_1,\ldots, E_m$ be the list of all such elliptic curves; this is a finite list, and for $N\le 707$ it can be easily extracted from the LMFDB.  If $J_H$ is completely decomposable there are nonnegative integers $e_1,\ldots,e_m$ for which
\begin{equation}\label{eq:LE}
L(J_H,s)= \prod_{i=1}^m L(E,s)^{e_i}.
\end{equation}
We may now assemble a list of $m$ primes $p_1,\ldots,p_m$ for which the $m\times m$ integer matrix $A:=[a_{p_i}(E_j)]_{ij}$ has nonzero determinant.  Such a list necessarily exists (the newforms corresponding to the $E_j$ are linearly independent), and in practice it does not take long to find. We now invert the matrix $A$ (over $\Q$), which we note depends only on the level~$N$, and can be computed before enumerating any subgroups.  As we enumerate subgroups $H\le \GL_2(\Z/N \Z)$ of level $N$ with $\det(H)=\widehat\Z^\times$ containing $-I$, for each $H$ we compute the integer vector $v:= (a_{p_1}(J_H),\ldots a_{p_m}(J_H))$ using the algorithm in \cite[\S 5]{RSZB} and then consider $A^{-1} v$.  If this vector does not lie in $\Z_{\ge 0}^m$ then $J_H$ cannot be completely decomposable, since its $L$-function cannot satisfy \eqref{eq:LE}.
For $N< 240$ the number of $E/\Q$ we need to consider is less than 1000, and the primes $p\le 8192$ are sufficient to obtain a matrix $[a_{p_i}(E_j)]$ of full rank from which we can extract a matrix with nonzero determinant by picking a subset of linearly independent rows.

If $A^{-1}v$ is a vector of nonnegative integers, this does not necessarily imply that $J_H$ is completely decomposable (although in practice this is almost always true).  If $J_H$ is not completely decomposable there certainly exists a list of primes $p_1,\ldots,p_m$ we could use to prove this, but it might not be the list we chose, and we don't want to take the time to check all possibilities up to the Sturm bound for $S_2(\Gamma_1(N) \cap \Gamma_0(N^2))$.  But while searching for unrealized genera we are happy to defer a rigorous verification to a later stage, since most of the $X_H$ we encounter will have genera that have already been realized (and when we find one that has not been realized, there will often be many, and we only need one example).

We implemented this strategy in \textsc{Magma}, using code from \cite{RSZB}; see the file \texttt{gl2split.m} in the GitHub repository linked to below.  We ran this algorithm on every level $N< 240$ (we included $N\le 70$ as a sanity check).  These computations were performed on a 256-core AMD EPYC 9754 2.25GHz CPU with 1.5TB RAM running Ubuntu 24.04 and took only a few days (less than a CPU-year of computation). This yielded a list of more than half a million modular Jacobians $J_H$ that appear to be completely decomposable. For each genus $g$ realized by one of these examples we determined the minimal~$H$ with $g(X_H)=g$, with the $H$ ordered according to their \href{https://beta.lmfdb.org/ModularCurve/Labels}{LMFDB label} $N.i.g.c.n$, whose first three parts encode the level $N$, index $i$, and genus~$g$ of $H$. These minimal examples are listed in Table~\ref{T:results}, and the full dataset is available at
\vspace{2pt}

\begin{center}
\footnotesize
\url{https://github.com/AndrewVSutherland/CompletelyDecomposableModularJacobians/}.
\end{center}
\vspace{2pt}

All but five of the minimal examples we found have level $N\le 70$, meaning the decomposition of $J_H$ is already available in the LMFDB and no further verification was required.  The exceptions all have level $120$, and after computing all the newforms in $S_2(\Gamma_1(120)\cap\Gamma_0(120^2))$, we were able to use the algorithm in \cite[\S 6]{RSZB} to rigorously determine the decomposition of $\Jac(X_H)$.  While $\dim(S_2(\Gamma_1(120)\cap\Gamma_0(120^2)) = 12032$ is rather large, it is the direct product of 279 newspaces $S_2(M,\chi)$ with $M|N^2$ of dimension at most $360$, so this was a feasible computation (the fact that $N=120$ is highly composite helps us here). The examples at level 120 yielded two previously unknown genera: 113 and 169.

\begin{center}
\footnotesize
\begin{table}[thp]
\setlength{\tabcolsep}{3.5pt}
\begin{tabular}{rlrlrlrl}
$g$ & LMFDB label & $g$ & LMFDB label & $g$ & LMFDB label & $g$ & LMFDB label\\\midrule
\textit{0} & \href{https://beta.lmfdb.org/ModularCurve/Q/1.1.0.a.1}{\texttt{1.1.0.a.1}} & \textit{21} & \href{https://beta.lmfdb.org/ModularCurve/Q/20.360.21.a.1}{\texttt{20.360.21.a.1}} & 45 & \href{https://beta.lmfdb.org/ModularCurve/Q/60.640.45.a.1}{\texttt{60.640.45.a.1}} & 85 & \href{https://beta.lmfdb.org/ModularCurve/Q/20.1440.85.b.1}{\texttt{20.1440.85.b.1}}\\
\textit{1} & \href{https://beta.lmfdb.org/ModularCurve/Q/6.6.1.a.1}{\texttt{6.6.1.a.1}} & 22 & \href{https://beta.lmfdb.org/ModularCurve/Q/20.360.22.a.1}{\texttt{20.360.22.a.1}} & \textit{49} & \href{https://beta.lmfdb.org/ModularCurve/Q/30.720.49.a.1}{\texttt{30.720.49.a.1}} & 89 & \texttt{120.1152.89.sn.1}\\
\textit{2} & \href{https://beta.lmfdb.org/ModularCurve/Q/10.30.2.a.1}{\texttt{10.30.2.a.1}} & 23 & \href{https://beta.lmfdb.org/ModularCurve/Q/30.360.23.a.1}{\texttt{30.360.23.a.1}} & 50 & \texttt{120.640.50.a.1} & \textit{97} & \href{https://beta.lmfdb.org/ModularCurve/Q/60.1440.97.y.1}{\texttt{60.1440.97.y.1}}\\
\textit{3} & \href{https://beta.lmfdb.org/ModularCurve/Q/7.168.3.a.1}{\texttt{7.168.3.a.1}} & 24 & \href{https://beta.lmfdb.org/ModularCurve/Q/24.384.24.a.1}{\texttt{24.384.24.a.1}} & 51 & \href{https://beta.lmfdb.org/ModularCurve/Q/60.720.51.u.1}{\texttt{60.720.51.u.1}} & 103 & \href{https://beta.lmfdb.org/ModularCurve/Q/60.1440.103.a.1}{\texttt{60.1440.103.a.1}}\\
\textit{4} & \href{https://beta.lmfdb.org/ModularCurve/Q/10.90.4.a.1}{\texttt{10.90.4.a.1}} & \textit{25} & \href{https://beta.lmfdb.org/ModularCurve/Q/24.384.25.a.1}{\texttt{24.384.25.a.1}} & 52 & \href{https://beta.lmfdb.org/ModularCurve/Q/60.720.52.i.1}{\texttt{60.720.52.i.1}} & 105 & \href{https://beta.lmfdb.org/ModularCurve/Q/60.1440.105.g.1}{\texttt{60.1440.105.g.1}}\\
\textit{5} & \href{https://beta.lmfdb.org/ModularCurve/Q/10.120.5.a.1}{\texttt{10.120.5.a.1}} & 26 & \href{https://beta.lmfdb.org/ModularCurve/Q/11.660.26.c.1}{\texttt{11.660.26.c.1}} & \textit{53} & \href{https://beta.lmfdb.org/ModularCurve/Q/60.720.53.a.1}{\texttt{60.720.53.a.1}} & 109 & \href{https://beta.lmfdb.org/ModularCurve/Q/60.1440.109.a.1}{\texttt{60.1440.109.a.1}}\\
\textit{6} & \href{https://beta.lmfdb.org/ModularCurve/Q/10.180.6.b.1}{\texttt{10.180.6.b.1}} & 27 & \href{https://beta.lmfdb.org/ModularCurve/Q/20.480.27.a.1}{\texttt{20.480.27.a.1}} & \textit{55} & \href{https://beta.lmfdb.org/ModularCurve/Q/36.864.55.a.1}{\texttt{36.864.55.a.1}} & \textbf{113} & \texttt{120.1440.113.hflr.1}\\
\textit{7} & \href{https://beta.lmfdb.org/ModularCurve/Q/10.180.7.a.1}{\texttt{10.180.7.a.1}} & 28 & \href{https://beta.lmfdb.org/ModularCurve/Q/18.486.28.a.1}{\texttt{18.486.28.a.1}} & \textit{57} & \href{https://beta.lmfdb.org/ModularCurve/Q/60.864.57.a.1}{\texttt{60.864.57.a.1}} & \textit{121} & \href{https://beta.lmfdb.org/ModularCurve/Q/60.1728.121.a.1}{\texttt{60.1728.121.a.1}}\\
\textit{8} & \href{https://beta.lmfdb.org/ModularCurve/Q/11.220.8.a.1}{\texttt{11.220.8.a.1}} & \textit{29} & \href{https://beta.lmfdb.org/ModularCurve/Q/20.480.29.a.1}{\texttt{20.480.29.a.1}} & \textit{61} & \href{https://beta.lmfdb.org/ModularCurve/Q/60.864.61.c.1}{\texttt{60.864.61.c.1}} & 129 & \texttt{120.1728.129.vwth.1}\\
\textit{9} & \href{https://beta.lmfdb.org/ModularCurve/Q/15.240.9.c.1}{\texttt{15.240.9.c.1}} & 31 & \href{https://beta.lmfdb.org/ModularCurve/Q/28.576.31.b.1}{\texttt{28.576.31.b.1}} & 65 & \href{https://beta.lmfdb.org/ModularCurve/Q/60.960.65.m.1}{\texttt{60.960.65.m.1}} & \textbf{135} & \href{https://beta.lmfdb.org/ModularCurve/Q/60.1920.135.a.1}{\texttt{60.1920.135.a.1}}\\
\textit{10} & \href{https://beta.lmfdb.org/ModularCurve/Q/15.180.10.a.1}{\texttt{15.180.10.a.1}} & 32 & \href{https://beta.lmfdb.org/ModularCurve/Q/30.540.32.a.1}{\texttt{30.540.32.a.1}} & \textbf{68} & \href{https://beta.lmfdb.org/ModularCurve/Q/60.960.68.a.1}{\texttt{60.960.68.a.1}} & \textbf{137} & \href{https://beta.lmfdb.org/ModularCurve/Q/60.1920.137.a.1}{\texttt{60.1920.137.a.1}}\\
\textit{11} & \href{https://beta.lmfdb.org/ModularCurve/Q/24.192.11.i.1}{\texttt{24.192.11.i.1}} & \textit{33} & \href{https://beta.lmfdb.org/ModularCurve/Q/24.576.33.a.1}{\texttt{24.576.33.a.1}} & 69 & \href{https://beta.lmfdb.org/ModularCurve/Q/60.960.69.c.1}{\texttt{60.960.69.c.1}} & \textbf{157} & \href{https://beta.lmfdb.org/ModularCurve/Q/60.2160.157.a.1}{\texttt{60.2160.157.a.1}}\\
12 & \href{https://beta.lmfdb.org/ModularCurve/Q/11.330.12.a.1}{\texttt{11.330.12.a.1}} & 34 & \href{https://beta.lmfdb.org/ModularCurve/Q/30.540.34.a.1}{\texttt{30.540.34.a.1}} & 72 & \href{https://beta.lmfdb.org/ModularCurve/Q/60.960.72.a.1}{\texttt{60.960.72.a.1}} & \textbf{159} & \href{https://beta.lmfdb.org/ModularCurve/Q/60.2160.159.gu.1}{\texttt{60.2160.159.gu.1}}\\
\textit{13} & \href{https://beta.lmfdb.org/ModularCurve/Q/10.360.13.b.1}{\texttt{10.360.13.b.1}} & 35 & \href{https://beta.lmfdb.org/ModularCurve/Q/60.480.35.a.1}{\texttt{60.480.35.a.1}} & \textit{73} & \href{https://beta.lmfdb.org/ModularCurve/Q/24.1152.73.a.1}{\texttt{24.1152.73.a.1}} & \textit{161} & \href{https://beta.lmfdb.org/ModularCurve/Q/48.2304.161.cax.1}{\texttt{48.2304.161.cax.1}}\\
14 & \href{https://beta.lmfdb.org/ModularCurve/Q/30.240.14.a.1}{\texttt{30.240.14.a.1}} & 36 & \href{https://beta.lmfdb.org/ModularCurve/Q/30.540.36.a.1}{\texttt{30.540.36.a.1}} & \textbf{75} & \href{https://beta.lmfdb.org/ModularCurve/Q/60.1080.75.gu.1}{\texttt{60.1080.75.gu.1}} & 163 & \href{https://beta.lmfdb.org/ModularCurve/Q/60.2160.163.c.1}{\texttt{60.2160.163.c.1}}\\
15 & \href{https://beta.lmfdb.org/ModularCurve/Q/20.240.15.a.1}{\texttt{20.240.15.a.1}} & \textit{37} & \href{https://beta.lmfdb.org/ModularCurve/Q/20.720.37.i.1}{\texttt{20.720.37.i.1}} & \textbf{76} & \href{https://beta.lmfdb.org/ModularCurve/Q/60.1080.76.dc.1}{\texttt{60.1080.76.dc.1}} & \textbf{169} & \texttt{120.2160.169.dcty.1}\\
16 & \href{https://beta.lmfdb.org/ModularCurve/Q/28.288.16.q.1}{\texttt{28.288.16.q.1}} & \textbf{38} & \href{https://beta.lmfdb.org/ModularCurve/Q/60.540.38.y.1}{\texttt{60.540.38.y.1}} & \textbf{77} & \href{https://beta.lmfdb.org/ModularCurve/Q/60.1080.77.ge.1}{\texttt{60.1080.77.ge.1}} & 193 & \href{https://beta.lmfdb.org/ModularCurve/Q/60.2880.193.dw.1}{\texttt{60.2880.193.dw.1}}\\
\textit{17} & \href{https://beta.lmfdb.org/ModularCurve/Q/20.360.17.e.1}{\texttt{20.360.17.e.1}} & 39 & \href{https://beta.lmfdb.org/ModularCurve/Q/60.540.39.e.1}{\texttt{60.540.39.e.1}} & \textbf{78} & \href{https://beta.lmfdb.org/ModularCurve/Q/60.1080.78.a.1}{\texttt{60.1080.78.a.1}} & \textit{205} & \href{https://beta.lmfdb.org/ModularCurve/Q/60.2880.205.a.1}{\texttt{60.2880.205.a.1}}\\
18 & \href{https://beta.lmfdb.org/ModularCurve/Q/60.240.18.m.1}{\texttt{60.240.18.m.1}} & 40 & \href{https://beta.lmfdb.org/ModularCurve/Q/60.540.40.a.1}{\texttt{60.540.40.a.1}} & 79 & \href{https://beta.lmfdb.org/ModularCurve/Q/60.1080.79.bg.1}{\texttt{60.1080.79.bg.1}} & 217 & \href{https://beta.lmfdb.org/ModularCurve/Q/60.2880.217.c.1}{\texttt{60.2880.217.c.1}}\\
\textit{19} & \href{https://beta.lmfdb.org/ModularCurve/Q/18.324.19.c.1}{\texttt{18.324.19.c.1}} & 41 & \href{https://beta.lmfdb.org/ModularCurve/Q/20.720.41.i.1}{\texttt{20.720.41.i.1}} & 81 & \href{https://beta.lmfdb.org/ModularCurve/Q/48.1152.81.id.1}{\texttt{48.1152.81.id.1}} & 325 & \href{https://beta.lmfdb.org/ModularCurve/Q/60.4320.325.a.1}{\texttt{60.4320.325.a.1}}\\
20 & \href{https://beta.lmfdb.org/ModularCurve/Q/20.360.20.a.1}{\texttt{20.360.20.a.1}} & \textit{43} & \href{https://beta.lmfdb.org/ModularCurve/Q/20.720.43.a.1}{\texttt{20.720.43.a.1}} & 82 & \href{https://beta.lmfdb.org/ModularCurve/Q/60.1080.82.a.1}{\texttt{60.1080.82.a.1}} & \textbf{409} & \href{https://beta.lmfdb.org/ModularCurve/Q/60.5760.409.c.1}{\texttt{60.5760.409.c.1}}\\
\bottomrule
\end{tabular}
\bigskip
\caption{Minimal modular curves of level $N<240$ with completely decomposable Jacobians ordered by genus $g$.  New genera are in bold, previously known modular genera are in italic. Unlinked labels are not yet in the LMFDB.}\label{T:results}
\end{table}
\end{center}

\begin{remark}
We believe that Table~\ref{T:results} includes every genus for which a completely decomposable modular Jacobian is known.  It includes every genus arising for $X_H$ of level $N<240$, and every genus arising from $X_0(N)$ for some $N$, the largest of which is $N=1200$ (we realize $g(X_0(1200))=205$ as $g(X_H)$ with $H$ of level $60$).  Table~\ref{T:results} does not include 32 genera for which nonmodular completely decomposable Jacobians are known, but not all of these~32 are realized over $\Q$.  This suggests two refinements of the question posed by Ekedahl and Serre: (1) in which genera are there completely decomposable modular Jacobians, and (2) in which genera are there nice curves $X/\Q$ with completely decomposable Jacobians?  Question (1) is a subset of (2), if we restrict to open $H\le \GL_2(\widehat\Z)$ with $\det(H)=\widehat\Z^\times$ as we do here.
\end{remark} 

\begin{remark}
Ekedahl and Serre also consider quotients of modular curves.  This allows them to realize several genera not known to occur for completely decomposable modular Jacobians; the genus 47 quotient $X_0(600)/\langle w_{24}\rangle$ is one such example.  We have not attempted to construct quotients of the completely decomposable modular Jacobians we found, but our search yielded more than 100,000 examples of genus $g>100$ whose quotients might yield new examples.  One could also consider modular Jacobians that are nearly but not completely decomposable, as these may admit quotients whose Jacobians are completely decomposable.
\end{remark}

\begin{remark}
The modular abelian variety $A_f$ associated to the Galois orbit of a non-CM newform is geometrically simple, so in most cases the $\Q$-isogeny decomposition of $J_H$ implied by \eqref{eq:LA} is a $\overline{\Q}$-isogeny decomposition (and in fact a decomposition over $\mathbb C$).  But the presence of complex multiplication complicates matters.  In our search we only checked for $J_H$ that are completely decomposable over $\Q$ and might have missed examples that are completely decomposable over $\overline{\Q}$.
\end{remark}

\bibliographystyle{alpha} 
\bibliography{references}

\newcommand{\etalchar}[1]{$^{#1}$}
\begin{thebibliography}{{LMF}25}

\bibitem[BBB{\etalchar{+}}21]{computing_cmf}
Alex~J. Best, Jonathan Bober, Andrew~R. Booker, Edgar Costa, John~E. Cremona,
  Maarten Derickx, Min Lee, David Lowry-Duda, David Roe, Andrew~V. Sutherland,
  and John Voight.
\newblock Computing classical modular forms.
\newblock In {\em Arithmetic geometry, number theory, and computation}, Simons
  Symp., pages 131--213. Springer, Cham, 2021.

\bibitem[BCFSe25]{magma}
W.~Bosma, J.~J. Cannon, C.~Fieker, and A.~Steel~(eds.).
\newblock Handbook of magma functions, version 2.28-17.
\newblock \url{https://magma.maths.usyd.edu.au/magma/handbook/}, 2025.

\bibitem[DS07]{char-p-3}
Claus Diem and Jasper Scholten.
\newblock Ordinary elliptic curves of high rank over {$\overline{\mathbb
  F}_p(x)$} with constant {$j$}-invariant. {II}.
\newblock {\em J. Number Theory}, 124(1):31--41, 2007.

\bibitem[EHR14]{EHR}
Noam~D. Elkies, Everett~W. Howe, and Christophe Ritzenthaler.
\newblock Genus bounds for curves with fixed {F}robenius eigenvalues.
\newblock {\em Proc. Amer. Math. Soc.}, 142(1):71--84, 2014.

\bibitem[ES93]{ekedahl-serre}
Torsten Ekedahl and Jean-Pierre Serre.
\newblock Exemples de courbes alg\'ebriques \`a{} jacobienne compl\`etement
  d\'ecomposable.
\newblock {\em C. R. Acad. Sci. Paris S\'er. I Math.}, 317(5):509--513, 1993.

\bibitem[Gon99]{char-p-2}
Josep Gonz\'alez.
\newblock Fermat {J}acobians of prime degree over finite fields.
\newblock {\em Canad. Math. Bull.}, 42(1):78--86, 1999.

\bibitem[{LMF}25]{lmfdb}
The {LMFDB Collaboration}.
\newblock The {L}-functions and modular forms database.
\newblock \url{https://beta.lmfdb.org}, 2025.
\newblock [Online; accessed 24 January 2025].

\bibitem[PR17]{paulhus-rojas}
Jennifer Paulhus and Anita~M. Rojas.
\newblock Completely decomposable {J}acobian varieties in new genera.
\newblock {\em Exp. Math.}, 26(4):430--445, 2017.

\bibitem[RR24]{anitarubi}
Rubí~E. Rodríguez and Anita~M. Rojas.
\newblock Decomposing abelian varieties into simple factors: algorithms and
  applications, 2024.

\bibitem[RSZB22]{RSZB}
Jeremy Rouse, Andrew~V. Sutherland, and David Zureick-Brown.
\newblock $\ell $ -adic images of {G}alois for elliptic curves over $\mathbb
  {Q}$ (and an appendix with {J}ohn {V}oight).
\newblock {\em Forum of Mathematics, Sigma}, 10:e62, 2022.

\bibitem[Ser20]{serre-book}
Jean-Pierre Serre.
\newblock {\em Rational points on curves over finite fields}, volume~18 of {\em
  Documents Math\'ematiques (Paris) [Mathematical Documents (Paris)]}.
\newblock Soci\'et\'e{} Math\'ematique de France, Paris, [2020] \copyright
  2020.
\newblock With contributions by Everett Howe, Joseph Oesterl\'e{} and
  Christophe Ritzenthaler.

\bibitem[Yam07]{yamauchi}
Takuya Yamauchi.
\newblock On {$\mathbb Q$}-simple factors of {J}acobian varieties of modular
  curves.
\newblock {\em Yokohama Math. J.}, 53(2):149--160, 2007.

\bibitem[Yui80]{char-p-1}
Noriko Yui.
\newblock On the {J}acobian variety of the {F}ermat curve.
\newblock {\em J. Algebra}, 65(1):1--35, 1980.

\end{thebibliography}

\end{document}